\documentclass[11pt,reqno]{amsart}

\usepackage{amsmath,amssymb,mathtools,esint}
\usepackage[letterpaper,margin=1.04in]{geometry}
\usepackage[hidelinks,pdfauthor={},pdfusetitle]{hyperref}
\newtheorem{theorem}{Theorem}[section]
\newtheorem{lemma}{Lemma}[section]
\newtheorem{corollary}{Corollary}[section]
\newtheorem{proposition}{Proposition}[section]
\newtheorem{remark}{Remark}[section]
\numberwithin{equation}{section}
\newcommand{\R}{\mathbb R}

\DeclareMathOperator{\Capacity}{Cap}
\newcommand{\W}{\mathcal W_k}
\newcommand{\V}{\mathcal V_k}
\allowdisplaybreaks[2]
\title{Hessian type equations with gradient source term and measure data}
\author{Shiguang Ma}
\address{School of Mathematical Science and LPMC,
Nankai University, Tianjin 300071,
People's Republic of China}
\email{msgdyx8741@nankai.edu.cn}

\author{Zijian Wang}
\address{Chern Institute of Mathematics and LPMC,
Nankai University, Tianjin 300071,
People's Republic of China}
\email{wzj@mail.nankai.edu.cn}

\begin{document}
\begin{abstract}
In this paper, we establish the necessary and sufficient conditions for the existence of nonnegative solutions to $F_k[-u]=|Du|^q\,dx+t\mu$ for some $t>0$, where $\mu$ is a finite nonnegative measure on $\mathbb{R}^N$. We consider the whole space $\mathbb{R}^N$ with $2\leq k<\frac{N}{2}$ and $q>\frac{Nk}{N-k}$, and a bounded domain $\Omega$ with $2\leq k<N$ and $q>k$. We also establish pointwise gradient estimates in terms of the Wolff potential.
\end{abstract}
\maketitle
\enlargethispage{6pt}
\tableofcontents
\section{Introduction}

In this paper, we establish necessary and sufficient conditions
for the existence of nonnegative solutions to the $k$-Hessian equation
with a gradient source and measure data
\begin{equation}\label{eq:problem}
 F_k[-u]=|Du|^q\,dx+t\mu
 \quad\text{in }\Omega,\quad u\ge0,
\end{equation}
where $2\le k<N$, $q>k$, $t>0$, $\Omega\subset\R^N$
is an open set, and $\mu$ is a nonnegative Radon measure
on $\Omega$. For $u\in C^2(\Omega)$, the classical $k$-Hessian operator
is defined by
\begin{equation*}
 F_k[-u]=\sigma_k(-D^2u)
 =\sum_{1\le i_1<\cdots<i_k\le N}
   \lambda_{i_1}\cdots\lambda_{i_k},
\end{equation*}
where $\lambda_1,\ldots,\lambda_N$ are the eigenvalues
of the Hessian matrix $-D^2u$. For $k$-Hessian equations, we refer to
\cite{CNS,GT,GU,IV,KR,TR,TW1,TW,TW3,TW4,UR,LA} for related results.

For the $k$-Hessian equation with a power source $u^q$,
Phuc and Verbitsky \cite{PV} established necessary and
sufficient conditions for the existence of solutions
in terms of nonlinear potentials and capacities.
For exponential reaction terms, Nguyen and V\'eron
\cite{NV} obtained existence results for quasilinear
and $k$-Hessian equations with measure data.
Ma and Wang \cite{MaWangPA} further studied the
$N$-Laplacian and $N/2$-Hessian equations with
exponential reaction terms and measure data.
 Liu, Ma and Zang \cite{LMZ} established
a pointwise gradient estimate for $k$-Hessian equations. For the quasilinear equation $-\Delta_p u=|Du|^q+\mu$, related solvability results were obtained by Phuc \cite{Phuc2010} and
Nguyen and Phuc \cite{NP}. And related quasilinear equations with source reaction terms
involving the function, its gradient and measure data is studied in \cite{Mixed}. 

Nonlinear potential theory plays an important role in the study of these equations. We first recall some notation needed to state our main results. For $\lambda=(\lambda_1,\ldots,\lambda_N)
\in\R^N$ and $1\le j\le k$, let
\begin{equation*}
 \sigma_j(\lambda)
 =\sum_{1\le i_1<\cdots<i_j\le N}
   \lambda_{i_1}\cdots\lambda_{i_j},
\end{equation*}
and define the G{\aa}rding cone
\begin{equation*}
 \Gamma_k=\{\lambda\in\R^N:
 \sigma_j(\lambda)>0,\ 1\le j\le k\}.
\end{equation*}

Following Trudinger and Wang \cite[Section~1]{TW}, we say that an upper semi-continuous function, $u:\Omega\to[-\infty,\infty)$,
is called $k$-convex in $\Omega$ if $F_k[P]\ge 0$ for all quadratic polynomials $P$ for which the difference $u-P$ has a finite local maximum
in $\Omega$. We will also call a $k$-convex function proper if it does not assume the value $-\infty$ identically on any component of
$\Omega$ and denote the class of proper $k$-convex functions in $\Omega$
by $\Phi^k(\Omega)$.

Trudinger and Wang first introduced the Hessian measure
for continuous $k$-convex functions in \cite{TW1} and
extended it to general functions in $\Phi_k(\Omega)$
in \cite[Theorem~1.1]{TW}.
For every $w\in\Phi_k(\Omega)$, the associated
$k$-Hessian measure $F_k[w]$ is uniquely determined by
\begin{equation*}
 \sigma_k(D^2w_j)\,dx\rightharpoonup F_k[w]
\end{equation*}
for any sequence of smooth $k$-convex functions
$w_j\to w$ in $L^1_{\mathrm{loc}}(\Omega)$.
Moreover, if $w\in C^2(\Omega)$, we have
\begin{equation*}
 dF_k[w]=\sigma_k(D^2w)\,dx.
\end{equation*}

Throughout this paper, we call a nonnegative function $u\in W^{1,q}_{\rm loc}(\Omega)$ is a solution of \eqref{eq:problem} if
$-u\in\Phi_k(\Omega)$ and
\begin{equation*}
 \int_\Omega\varphi\,dF_k[-u]
 =\int_\Omega\varphi|Du|^q\,dx
  +t\int_\Omega\varphi\,d\mu
 \quad\text{for every }\varphi\in C_c(\Omega).
\end{equation*}
Here $Du$ is the weak gradient of $u$,
$C_c(\Omega)$ denotes the space of continuous
compactly supported functions on $\Omega$,
and $dx$ denotes the Lebesgue measure.

We next introduce the Riesz and Bessel capacities used
in our solvability conditions; see \cite{AH,PV}.
Let $I_k$ be the Riesz kernel of order $k$ and
$G_k$ be the Bessel kernel of the same order, which is defined by
\begin{equation}
 I_k(x)=c_{N,k}|x|^{k-N},
 \quad
 \widehat G_k(\xi)=(1+|\xi|^2)^{-k/2},
\end{equation}
where $c_{N,k}>0$ is the usual normalization constant and $\widehat G_k$ denotes the Fourier transform of $G_k$. For $J=I_k$ or $G_k$, $s>1$, and a compact set $K\subset\R^N$, the corresponding capacity is
\begin{equation}\label{eq:cap-def}
 \Capacity_{J,s}(K)
 =\inf\left\{
 \|f\|_{L^s(\R^N)}^s:
 f\in L^s(\R^N),\ f\ge0,\
 J*f\ge1\text{ on }K
 \right\}.
\end{equation}

Finally, for a nonnegative Radon measure $\lambda$
on $\R^N$ and $0<L\le\infty$, we define
the truncated potentials
\begin{equation}\label{eq:wdef}
\begin{aligned}
 \W^L[\lambda](x)
 &=\int_0^L
 \left(\frac{\lambda(B_r(x))}{r^{N-k}}\right)^{1/k}
 \frac{dr}{r},
 &\W&=\W^\infty,\\
 \V^L[\lambda](x)
 &=\int_0^L
 \left(\frac{\lambda(B_r(x))}{r^{N-2k}}\right)^{1/k}
 \frac{dr}{r},
 &\V&=\V^\infty.
\end{aligned}
\end{equation}
Here $B_r(x)$ is the Euclidean ball of radius $r$
centered at $x$. These correspond to the Wolff
potentials with parameters
$(k/(k+1),k+1)$ and $(2k/(k+1),k+1)$, respectively;
see \cite[Section~2]{PV}.
Measures defined on $\Omega$ are extended by zero
outside $\Omega$ whenever these potentials are used.

\begin{theorem}\label{thm:exist}
Let $2\le k<N/2$, $q>Nk/(N-k)$, and let $\mu$ be a nonnegative
finite measure on $\R^N$. Then the following statements are equivalent.
\begin{enumerate}
\item For some $t>0$, \eqref{eq:problem} has a solution in $\R^N$
with $\inf_{\R^N}u=0$.
\item There exists a constant $D<\infty$ such that for any compact set $K\subset\R^N$
\begin{equation}\label{eq:riesz-condition}
 \mu(K)\le D\Capacity_{I_k,q/(q-k)}(K).
\end{equation}
\end{enumerate}
If $D\le\varepsilon_0(N,k,q)$, there exists a solution for $t=1$ with
\begin{equation}\label{eq:global-gradient}
 |Du|\le C\W[\mu]\quad\text{a.e. in }\R^N,\quad Du\in L^q(\R^N).
\end{equation}
Conversely, every solution for $t=1$ satisfies that for any compact set $K\subset\R^N$,
\begin{equation}
 F_k[-u](K)\le C(N,k,q)\Capacity_{I_k,q/(q-k)}(K).  
\end{equation}
\end{theorem}
In fact, we can extend the existence criterion of Theorem~\ref{thm:exist}
to locally finite measures and obtain a two-sided pointwise estimate
for $u$ in terms of the Wolff potential of $\mu$.
\begin{corollary}\label{cor:locally-finite}
Suppose that $2\le k<N/2$, $q>Nk/(N-k)$, and $\mu$ is a nonnegative
locally finite Radon measure on $\R^N$. Then, there exist a constant $t>0$ and a
nonnegative solution of $F_k[-u]=|Du|^q\,dx+t\mu$ with
$\inf_{\R^N}u=0$ if and only if \eqref{eq:riesz-condition} holds
with $D<\infty$ and
\begin{equation}\label{eq:tail-locally-finite}
 \int_1^\infty\left(\frac{\mu(B_r(0))}{r^{N-2k}}\right)^{1/k}
 \frac{dr}{r}<\infty .
\end{equation}
Moreover, if \eqref{eq:tail-locally-finite} and
\eqref{eq:riesz-condition} holds with $D\le\varepsilon_0(N,k,q)$,
then for $t=1$, \eqref{eq:problem} has a nonnegative solution $u$
with $\inf_{\R^N}u=0$ such that for every $x\in\R^N$, we have
\begin{equation*}
c\V[\mu](x)\le u(x)\le C\V[\mu](x),
\end{equation*}
and
\begin{equation*}
|Du(x)|\le C\W[\mu](x),
\end{equation*}
for almost every $x\in\R^N$. If $q<2k$, then \eqref{eq:riesz-condition} implies
\eqref{eq:tail-locally-finite}.
\end{corollary}

On a bounded domain $\Omega$ we consider the Dirichlet problem
\begin{equation}\label{eq:dirichlet}
\begin{cases}
 F_k[-u]=|Du|^q\,dx+t\mu & \text{in }\Omega,\\
 u=0 & \text{on }\partial\Omega.
\end{cases}
\end{equation}
The boundary condition means that $\lim_{\Omega\ni x\to z}u(x)=0$
for every $z\in\partial\Omega$.

\begin{theorem}\label{thm:ball}
Let $2\le k<N$, $q>k$, and let $\Omega$ be a smooth bounded uniformly
$(k-1)$-convex domain. Suppose that $\mu$ is a nonnegative finite measure
supported in a compact set $E\Subset\Omega$. Then, the following statements are
equivalent.
\begin{enumerate}
\item For some $t>0$, the Dirichlet problem \eqref{eq:dirichlet} has a
nonnegative solution.
\item There exists a constant $D<\infty$ such that for every compact set $K\subset E$
\begin{equation}\label{eq:bessel-condition}
 \mu(K)\le D\Capacity_{G_k,q/(q-k)}(K).
\end{equation}
\end{enumerate}
If $D\le\varepsilon_0(N,k,q,\Omega,E)$, there exists a solution for $t=1$ with
\begin{equation}\label{eq:bounded-gradient}
 u\in W_0^{1,q}(\Omega),\quad
 |Du|\le C\W^{8\operatorname{diam}(\Omega)}[\mu]\quad\text{a.e. in }\Omega.
\end{equation}
Conversely, every solution for $t=1$ satisfies, for every compact set $K\subset E$,
\begin{equation}
 F_k[-u](K)\le C(N,k,q,\Omega,E)\Capacity_{G_k,q/(q-k)}(K).
\end{equation}
\end{theorem}

\begin{remark}
For every nonnegative finite measure $\mu$ supported in $E$,
\eqref{eq:bessel-condition} always holds for some $D<\infty$
when $k<q<Nk/(N-k)$. In fact, for $s=q/(q-k)$,
\begin{equation*}
 \Capacity_{G_k,s}(\{x\})>0
 \quad\Longleftrightarrow\quad
 ks>N
 \quad\Longleftrightarrow\quad
 q<\frac{Nk}{N-k};
\end{equation*}
see \cite[Theorem~2.5.1]{AH}.
Moreover, Proposition~\ref{prop:hessian-capacity} shows
that if $q\ge Nk/(N-k)$ and $\mu(\{x_0\})>0$
for some $x_0\in\Omega$, then \eqref{eq:problem} admits no nonnegative solution for any $t>0$.
\end{remark}

\section{The Capacity Estimates}\label{sec:necessary}
In this section, we prove the following capacity estimates for the Hessian measure $F_k[-u]$ of any nonnegative solution $u$ of \eqref{eq:problem}. These estimates also imply that $\mu$ satisfies the capacity conditions in Theorems~\ref{thm:exist} and~\ref{thm:ball}. As an application of these capacity estimates, at the end of this Section, a Liouville-type theorem for entire supersolutions will be proved.
\begin{proposition}\label{prop:hessian-capacity}
Let $2\le k<N$, $q>k$, and $s=q/(q-k)$. Fix $t>0$.
\begin{enumerate}
\item Suppose that $ks<N$ and $u$ is a nonnegative solution of
\eqref{eq:problem} in $\R^N$. Then for every compact set
$K\subset\R^N$ we have
\begin{equation*}
t\mu(K)\le F_k[-u](K)
\le C(N,k,q)\Capacity_{I_k,s}(K).
\end{equation*}

\item Let $\Omega\subseteq\R^N$ be a open set, suppose that $u$
is a nonnegative solution of \eqref{eq:problem} in $\Omega$.
Then for every compact set $E\Subset\Omega$ and every compact set
$K\subset E$ we have
\begin{equation*}
t\mu(K)\le F_k[-u](K)
\le C(N,k,q,E,\Omega)\Capacity_{G_k,s}(K).
\end{equation*}
\end{enumerate}
\end{proposition}
\begin{lemma}\label{lem:hessian-cutoff}
Let $2\le k<N$, $q>k$, $s=q/(q-k)$, and let $m$ be an integer with
$m>ks$. Suppose that $w\in\Phi_k(\Omega)\cap W^{1,q}_{\mathrm{loc}}(\Omega)$
and $F_k[w]\ge |Dw|^q\,dx$ in $\Omega$. Then for every
$\varphi\in C_c^\infty(\Omega)$ with $0\le\varphi\le1$, we have
\begin{equation}\label{eq:hessian-cutoff}
 \int_\Omega\varphi^m\,dF_k[w]
 \le C(N,k,q,m)\int_\Omega |D^k\varphi|^s\,dx .
\end{equation}
\end{lemma}

\begin{proof}
We first assume that $w\in C^\infty(\Omega)\cap\Phi_k(\Omega)$. Define that $F_0=1$, $F_\ell=\sigma_\ell(D^2w)$ and
$F_\ell^{ij}=\partial\sigma_\ell(D^2w)/\partial w_{ij}$. Fix $1\le\ell\le k$ and take $p=2+(k-\ell)/\ell$. We have
\begin{equation*}
 |Dw|^{k-\ell}F_\ell
 \le \frac{N}{k}\,
 \sigma_\ell\bigl(D(|Dw|^{p-2}Dw)\bigr).
\end{equation*}
In fact, for $\ell<k$ notice that
$\ell(p-2)=k-\ell<N(k-\ell)/(N-k)$, so it follows from \cite[(4.12)]{TW}. For $\ell=k$, it reduces to $F_k\le(N/k)F_k$ so we get it immediately.

Let $0\le h\in C_c^\infty(\Omega)$. Multiplying the preceding inequality by $h$ and integrating over $\Omega$, we obtain
\begin{equation*}
\begin{aligned}
 \int_\Omega h|Dw|^{k-\ell}F_\ell\,dx
 \le \frac{N}{k}\int_\Omega
       h\,\sigma_\ell\bigl(D(|Dw|^{p-2}Dw)\bigr)\,dx=-\frac{N}{k\ell}\int_\Omega |Dw|^{k-\ell}
       \sum_{i,j=1}^N F_\ell^{ij}\,\partial_i h\,\partial_j w\,dx .
\end{aligned}
\end{equation*}
The last equality follows from the formula of Trudinger and Wang \cite[(4.17)]{TW}.

It remains to estimate the sum in the last integral. Since
$w\in\Phi_k(\Omega)$ and $\ell\le k$, the matrix $(F_\ell^{ij})$
is positive semidefinite. Its trace is
\begin{equation*}
 \sum_{i=1}^N F_\ell^{ii}
 =\left.\frac{d}{dt}\sigma_\ell(D^2w+tI)\right|_{t=0}
 =(N-\ell+1)F_{\ell-1}.
\end{equation*}
At each point, we can choose orthonormal coordinates such that $(F_\ell^{ij})$ is diagonal matrix, then since $F_\ell^{ii}\ge0$ for every $i$,
the Cauchy--Schwarz inequality gives
\begin{equation*}
\begin{aligned}
 \left|\sum_{i,j=1}^N F_\ell^{ij}\,\partial_i h\,\partial_j w\right|
 &\le \left(\max_{1\le i\le N} F_\ell^{ii}\right)
       \sum_{i=1}^N |\partial_i h|\,|\partial_i w|\\
 &\le \left(\sum_{i=1}^N F_\ell^{ii}\right)|Dh|\,|Dw|
  \\&=(N-\ell+1)F_{\ell-1}|Dh|\,|Dw|.
\end{aligned}
\end{equation*}
Then for such $h$, this estimate gives us
\begin{equation}\label{eq:weighted-hessian-step}
\begin{aligned}
 \int_\Omega h|Dw|^{k-\ell}F_\ell\,dx
 \le\frac{N(N-\ell+1)}{k\ell}
       \int_\Omega |Dh|\,|Dw|^{k-\ell+1}F_{\ell-1}\,dx .
\end{aligned}
\end{equation}
By approximation, \eqref{eq:weighted-hessian-step} also holds for every nonnegative Lipschitz function $h$ with compact support in $\Omega$.

We now iterate \eqref{eq:weighted-hessian-step} to reduce $F_k$
to $F_0=1$ and we use the convention that $D^0\varphi=\varphi$. We claim that there exists a constant $C=C(N,k,m)$, such that
for every integer $0\le d\le k$,
\begin{equation*}
 \int_\Omega \varphi^m F_k\,dx
 \le C
\sum_{\substack{\alpha_1,\ldots,\alpha_m\ge0\\\alpha_1+\cdots+\alpha_m=d}}
 \int_\Omega |Dw|^d F_{k-d}
       \prod_{i=1}^m |D^{\alpha_i}\varphi|\,dx ,
\end{equation*}
where the $\alpha_i$ are integers. 

For $d=0$, all $\alpha_i$ vanish
and the claim holds with equality. For $d=1$, taking
$h=\varphi^m$ and $\ell=k$ in \eqref{eq:weighted-hessian-step} we obtain
\begin{equation*}
 \int_\Omega \varphi^m F_k\,dx
 \le C \int_\Omega
       \varphi^{m-1}|D\varphi|\,|Dw|F_{k-1}\,dx .
\end{equation*}
Suppose the claim holds for some $1\le d<k$. For any nonnegative
integers $\alpha_1,\ldots,\alpha_m$ with $\alpha_1+\cdots+\alpha_m=d$,
the product $\prod_{i=1}^m |D^{\alpha_i}\varphi|$ is a nonnegative
Lipschitz function supported in $\Omega$.
By the product rule, almost everywhere in $\Omega$,
\begin{equation*}
 \left|D\left(\prod_{i=1}^m |D^{\alpha_i}\varphi|\right)\right|
 \le \sum_{j=1}^m |D^{\alpha_j+1}\varphi|
                  \prod_{i\ne j}|D^{\alpha_i}\varphi| .
\end{equation*}
Then, by the induction hypothesis, using \eqref{eq:weighted-hessian-step}
with $\ell=k-d$ and $h=\prod_{i=1}^m |D^{\alpha_i}\varphi|$,
we obtain
\begin{equation*}
\begin{aligned}
 \int_\Omega \varphi^m F_k\,dx
 &\quad\le C\sum_{\alpha_1+\cdots+\alpha_m=d}
       \int_\Omega |Dw|^d F_{k-d}
       \prod_{i=1}^m |D^{\alpha_i}\varphi|\,dx\\
 &\quad\le C\sum_{\alpha_1+\cdots+\alpha_m=d}
       \int_\Omega |Dw|^{d+1}F_{k-d-1}
       \left|D\left(\prod_{i=1}^m |D^{\alpha_i}\varphi|\right)\right|\,dx\\
 &\quad\le C\sum_{\alpha_1+\cdots+\alpha_m=d}\sum_{j=1}^m
       \int_\Omega |Dw|^{d+1}F_{k-d-1}
       |D^{\alpha_j+1}\varphi|
       \prod_{i\ne j}|D^{\alpha_i}\varphi|\,dx\\
 &\quad\le C\sum_{\beta_1+\cdots+\beta_m=d+1}
       \int_\Omega |Dw|^{d+1}F_{k-d-1}
       \prod_{i=1}^m |D^{\beta_i}\varphi|\,dx .
\end{aligned}
\end{equation*}
This proves the claim for $d+1$, then we prove the claim. Taking $d=k$ in the claim and using $F_0=1$, we obtain
\begin{equation*}
\begin{aligned}
 \int_\Omega \varphi^m\,dF_k[w]
 \le C\sum_{\alpha_1+\cdots+\alpha_m=k}
       \int_\Omega |Dw|^k
       \prod_{i=1}^m |D^{\alpha_i}\varphi|\,dx .
\end{aligned}
\end{equation*}

We next extend the estimate to
$w\in\Phi_k(\Omega)\cap W^{1,q}_{\rm loc}(\Omega)$. Fix an open set $U\Subset\Omega$ containing $\operatorname{supp}\varphi$,
and let $w_j$ be the standard mollifications of $w$. By
\cite[Lemma~2.3]{TW}, $w_j\in C^\infty(U)\cap\Phi_k(U)$ for all
sufficiently large $j$. Since $w\in W^{1,q}_{\rm loc}(\Omega)$,
we have $Dw_j\to Dw$ in $L^q(U)$. Applying the estimate for $d=k$ to $w_j$ and letting $j\to\infty$,
we obtain
\begin{equation}\label{eq:iterated-hessian-cutoff}
\begin{aligned}
 \int_\Omega \varphi^m\,dF_k[w]
 &=\lim_{j\to\infty}\int_U \varphi^m\,dF_k[w_j]\\
 &\le C\lim_{j\to\infty}
       \sum_{\alpha_1+\cdots+\alpha_m=k}
       \int_U |Dw_j|^k\prod_{i=1}^m |D^{\alpha_i}\varphi|\,dx\\
 &=C\sum_{\alpha_1+\cdots+\alpha_m=k}
       \int_\Omega |Dw|^k\prod_{i=1}^m |D^{\alpha_i}\varphi|\,dx .
\end{aligned}
\end{equation}
The first equality follows from \cite[Theorem~1.1(ii)]{TW} with
test function $\varphi^m$. For the last equality, $q>k$ gives $Dw_j\to Dw$ in $L^k(U)$ and hence $|Dw_j|^k\to|Dw|^k$ in $L^1(U)$.
Each product of derivatives of $\varphi$ is bounded and supported in $U$.

It remains to estimate the terms on the right hand side of \eqref{eq:iterated-hessian-cutoff}.
Fix nonnegative integers $\alpha_1,\ldots,\alpha_m$ with
$\alpha_1+\cdots+\alpha_m=k$. At least $m-k$ of the $\alpha_i$ are
zero. Since $0\le\varphi\le1$ and
$m-k-mk/q=m/s-k>0$, we have
\begin{equation*}
 \prod_{i=1}^m |D^{\alpha_i}\varphi|
 \le \varphi^{m-k}\prod_{\alpha_i>0}|D^{\alpha_i}\varphi|
 \le \varphi^{mk/q}\prod_{\alpha_i>0}|D^{\alpha_i}\varphi| .
\end{equation*}
Thus
\begin{equation*}
\begin{aligned}
 \int_\Omega |Dw|^k\prod_{i=1}^m|D^{\alpha_i}\varphi|\,dx
 &\le \int_\Omega (\varphi^m|Dw|^q)^{k/q}
                \prod_{\alpha_i>0}|D^{\alpha_i}\varphi|\,dx\\
 &\le \left(\int_\Omega\varphi^m|Dw|^q\,dx\right)^{k/q}
       \prod_{\alpha_i>0}
       \|D^{\alpha_i}\varphi\|_{L^{ks/\alpha_i}(\Omega)}\\
 &\le C\left(\int_\Omega\varphi^m|Dw|^q\,dx\right)^{k/q}
       \|D^k\varphi\|_{L^s(\Omega)} .
\end{aligned}
\end{equation*}
The second inequality follows from H\"older's inequality, since
$
 \frac{k}{q}+\sum_{\alpha_i>0}\frac{\alpha_i}{ks}
 =\frac{k}{q}+\frac1s=1 .
$
For the last inequality, the Gagliardo--Nirenberg inequality applied
to the zero extension of $\varphi$ gives, for $1\le j<k$,
\begin{equation*}
 \|D^j\varphi\|_{L^{ks/j}(\Omega)}
 \le C\|\varphi\|_{L^\infty(\Omega)}^{1-j/k}
       \|D^k\varphi\|_{L^s(\Omega)}^{j/k}
 \le C\|D^k\varphi\|_{L^s(\Omega)}^{j/k}.
\end{equation*}
For $j=k$ the last bound holds with $C=1$, and
$\sum_{\alpha_i>0}\alpha_i/k=1$.

Summing over $\alpha_1+\cdots+\alpha_m=k$ in
\eqref{eq:iterated-hessian-cutoff}, we obtain
\begin{equation*}
\begin{aligned}
 \int_\Omega\varphi^m\,dF_k[w]
 &\le C\left(\int_\Omega\varphi^m|Dw|^q\,dx\right)^{k/q}
       \|D^k\varphi\|_{L^s(\Omega)}\\
 &\le C\left(\int_\Omega\varphi^m\,dF_k[w]\right)^{k/q}
       \|D^k\varphi\|_{L^s(\Omega)}\\
 &\le \frac12\int_\Omega\varphi^m\,dF_k[w]
       +C\int_\Omega |D^k\varphi|^s\,dx .
\end{aligned}
\end{equation*}
In the second inequality we use $|Dw|^q\,dx\le F_k[w]$, and the last inequality follows from Young's inequality with exponents $q/k$ and $s$.
Subtracting the first term on the right, then we finish the proof.
\end{proof}
The following capacity estimate is proved in \cite{Mazya2011}
\begin{lemma}\label{lem:cutoff-cap}
Let $1\le k<N$ be an integer, $1<s<\infty$, and let $K\subset\R^N$ be a compact set.
Then, if $ks<N$, we have
\begin{equation}\label{eq:cutoff-cap-riesz}
 \inf_{\substack{\varphi\in C_c^\infty(\R^N),\ 0\le\varphi\le1\\
                   \varphi=1\ \mathrm{near}\ K}}
       \int_{\R^N}|D^k\varphi|^s\,dx
 \le C(N,k,s)\operatorname{Cap}_{I_k,s}(K).
\end{equation}
If $K\subset E\Subset\Omega$ for a fixed compact set $E$, then
\begin{equation}\label{eq:cutoff-cap-bessel}
 \inf_{\substack{\varphi\in C_c^\infty(\Omega),\ 0\le\varphi\le1\\
                   \varphi=1\ \mathrm{near}\ K}}
        \int_\Omega |D^k\varphi|^s\,dx
 \le C(N,k,s,E,\Omega)\operatorname{Cap}_{G_k,s}(K).
\end{equation}
\end{lemma}

We now apply the preceding capacity estimates to the Hessian
measure.

\begin{proof}[Proof of Proposition~\ref{prop:hessian-capacity}]
Set $w=-u$ and $m=\lfloor ks\rfloor+1$.
Since
$F_k[w]\ge |Dw|^q\,dx$,
Lemma~\ref{lem:hessian-cutoff} gives
\begin{equation*}
\begin{aligned}
 t\mu(K)
 \le F_k[w](K)
& \le \inf_{\substack{
 \varphi\in C_c^\infty(\Omega),\ 0\le\varphi\le1
 \varphi=1\ \mathrm{near}\ K}}
 \int_\Omega\varphi^m\,dF_k[w]
\\&\le C\inf_{\substack{
 \varphi\in C_c^\infty(\Omega),\ 0\le\varphi\le1\\
 \varphi=1\ \mathrm{near}\ K}}
 \int_\Omega|D^k\varphi|^s\,dx.
\end{aligned}
\end{equation*}
The conclusions follow from
\eqref{eq:cutoff-cap-riesz} and
\eqref{eq:cutoff-cap-bessel}.
\end{proof}
We conclude this section with a Liouville-type
theorem for entire supersolutions
\begin{corollary}\label{cor:liouville}
Let $2\le k<N$, $q>k$, and suppose that
$w\in\Phi_k(\R^N)\cap W^{1,q}_{\rm loc}(\R^N)$ satisfies
$F_k[w]\ge |Dw|^q\,dx$. Then $w$ is constant if $q\le Nk/(N-k)$.
For $2k\ge N$, the same conclusion holds if $q>k$ and $w\le0$.
\end{corollary}

\begin{proof}
Put $s=q/(q-k)$. If $q<Nk/(N-k)$, then $ks>N$. Fix a function $\chi\in C_c^\infty(\R^N)$
satisfies that $0\le\chi\le1$ in and $\chi=1 $ on $B_1$.
By \eqref{eq:hessian-cutoff} with $\varphi(x)=\chi(x/R)$,
\begin{equation*}
 F_k[w](B_R)\le C R^{N-ks}\int_{\R^N}|D^k\chi|^s\,dx.
\end{equation*}
Then the right-hand side tends to zero as $R\to\infty$. 

If $q=Nk/(N-k)$, then $ks=N$. Fix $a>0$ and let $R>4a$.
By \cite[\S13.1.2, Proposition~2, proof of (13.1.5)]{Mazya2011}, we may choose $\varphi\in C_c^\infty(B_R)$ with $0\le\varphi\le1$ and $\varphi=1$ on $B_a$ such that
\begin{equation*}
 \int_{B_R}|D^k\varphi|^s\,dx\le C(\log(R/a))^{1-s},
\end{equation*}
where $C$ is independent of $R$. Applying \eqref{eq:hessian-cutoff}
with this $\varphi$, we obtain
\begin{equation*}
 F_k[w](B_a)\le C(\log(R/a))^{1-s}.
\end{equation*}
Since $s>1$, letting $R\to\infty$ yields $F_k[w](B_a)=0$.

Suppose $2k\ge N$, $w\le0$ and $F_k[w]\ne0$. Choose $x$ such that $-w(x)<\infty$ and $r_0>0$ with $F_k[w](B_{r_0}(x))>0$.
The lower bound of Labutin \cite[Theorem~7.2]{PV} on $B_{24R}(x)$ gives
\begin{equation*}
 -w(x)\ge c\int_{r_0}^R
    \left(\frac{F_k[w](B_r(x))}{r^{N-2k}}\right)^{1/k}\frac{dr}{r}
 \ge c F_k[w](B_{r_0}(x))^{1/k}
             \int_{r_0}^R r^{2-N/k}\frac{dr}{r}.
\end{equation*}
Since $2k\ge N$, the last integral tends to infinite as $R\to\infty$,
contradicting to the choice of $x$ with $-w(x)<\infty$.
Hence $F_k[w]=0$.

In each case $|Dw|^q\,dx\le F_k[w]=0$, so the subharmonic representative
of $w$ is a constant.
\end{proof}

\section{Hessian Type Equation In $\mathbb{R}^N$}
\label{sec:existence}
We first recall the following estimates for the Wolff potential defined in \eqref{eq:wdef}.
\begin{lemma}\label{lem:potential}
Let $2\le k<N$, $q>k$, $s=q/(q-k)$, and let $\mu$ be a nonnegative
finite measure on $\R^N$. Suppose that one of the following holds.
\begin{enumerate}
\item $L=\infty$, $ks<N$, and
$\mu(K)\le D\Capacity_{I_k,s}(K)$ for every compact set
$K\subset\R^N$.
\item $0<L<\infty$,
$\operatorname{supp}\mu\subset\overline B_{L/8}(x_0)$,
and $\mu(K)\le D\Capacity_{G_k,s}(K)$ for every compact set
$K\subset\R^N$.
\end{enumerate}
Then
\begin{equation}\label{eq:energy}
 \int_{\R^N}\bigl(\W^L[\mu]\bigr)^q\,dx
 \le C D^{(q-k)/k}\mu(\R^N),
\end{equation}
and
\begin{equation}\label{eq:closure}
 \W^L\bigl[(\W^L[\mu])^q\,dx\bigr](x)
 \le C D^{(q-k)/k^2}\W^L[\mu](x).
\end{equation}
Estimate \eqref{eq:closure} holds almost everywhere in $\R^N$
in case \textup{(1)} and almost everywhere in $B_{L/4}(x_0)$
in case \textup{(2)}.
The constant $C$ depends only on $N,k,q$ and additionally on $L$ in case \textup{(2)}.
\end{lemma}
\begin{proof}
The proof follows from \cite[Propositions~2.4 and~2.5]{BVNV} immediately.
\end{proof}
We next recall the pointwise gradient estimate of Liu, Ma and Zang \cite[Theorem~1.2]{LMZ}, which will be used later.
\begin{proposition}\label{prop:gradient-estimate}
For $2\le k<N$, let $\Omega\subset\R^N$ be an
open set, and $w\in\Phi_k(\Omega)$.
Set $\nu=F_k[w]$. If $B_r(x)\Subset\Omega$
and $x$ is a Lebesgue point of $Dw$, then
\begin{equation}\label{eq:lmz}
 |Dw(x)|
 \le \frac{C}{r}\fint_{B_r(x)}|w(y)|\,dy
       +C\W^r[\nu](x),
\end{equation}
where $C=C(N,k)$.
\end{proposition}
Then we can begin to prove Theorem~\ref{thm:exist}
\begin{proof}[Proof of Theorem~\ref{thm:exist}]
We first prove existence theorem for $t=1$ when $D$ in \eqref{eq:riesz-condition}
is sufficiently small. If $\mu=0$, take $u=0$ then we finish the proof. Suppose $\mu\ne0$.
Since $q>Nk/(N-k)$, we have $kq/(q-k)<N$, and \eqref{eq:energy} gives
$\W[\mu]\in L^q(\R^N)$.

We apply the Schauder Fixed Point Theorem to prove the theorem. Write $B_m=B_m(0)$. Let $m\ge2$ be an integer with $\mu(B_m)>0$,  a constant $A>0$ is chosen independently of $m$,
we define that
\begin{equation*}
 \mathcal K_m=\{g\in L^q(B_m;\R^N):
                    |g|\le A\W[\mu]\ \text{a.e. in }B_m\}.
\end{equation*}
Then it is a nonempty closed convex subset of $L^q(B_m;\R^N)$.
Fix a nonnegative even mollifier $\rho\in C_c^\infty(B_1)$ with
$\int_{\R^N}\rho\,dx=1$, and let $\rho_h=h^{-N}\rho(\cdot/h)$.
For $g\in\mathcal K_m$, we define that
\begin{equation}\label{notation in thm1}
 \lambda_g=\mathbf1_{B_m}(\mu+|g|^q\,dx),\quad
 b_m=m^{-N}\mu(B_m),\quad f_g=\rho_{1/m}*\lambda_g+b_m.
\end{equation}
Then $f_g\in C^\infty(\overline B_{4m})$ and $f_g\ge b_m>0$. Therefore, by \cite[Theorem~3]{CNS}, the equation
\begin{equation*}
 \sigma_k(D^2w_g)=f_g\quad\text{in }B_{4m},\quad
 w_g=0\quad\text{on }\partial B_{4m}.
\end{equation*}
has a unique smooth admissible solution.

Set $T_mg=Dw_g|_{B_m}$ and $\sigma_g=\mathbf1_{B_{4m}}f_g\,dx$.
For $x\in B_m$, applying \eqref{eq:lmz} with radius $m$, we obtain
\begin{equation*}
 |T_mg(x)|=|Dw_g(x)|
 \le \frac{C}{m}\fint_{B_m(x)}|w_g(y)|\,dy
      +C\W^m[\sigma_g](x).
\end{equation*}
To prove $T_mg\in\mathcal K_m$, we need to bound the right-hand side
by $A\W[\mu](x)$ almost everywhere in $B_m$.

We first estimate the average of $|w_g|$. Since $w_g$ is a admissible solution, then $\Delta w_g\ge0$ in $B_{4m}$. Moreover, notice that $w_g=0$ on $\partial B_{4m}$, the maximum principle gives
$w_g\le0$ in $B_{4m}$. Applying \cite[Theorem~7.4]{PV} to
$u=-w_g$, with $F_k[-u]=\sigma_g$ in $B_{4m}$ and zero boundary
values, we obtain that for any $x\in B_{4m}$
\begin{equation} \label{wg V estimate}
 |w_g(x)|=-w_g(x)
 \le C\V^{16m}[\sigma_g](x)
 \le C\V[\sigma_g](x),
\end{equation}
where $16m=2\operatorname{diam}(B_{4m})$. For every nonnegative finite measure $\sigma$ on $\R^N$, $x\in\R^N$,
and $r,S>0$. By direct calculation, Tonelli's theorem gives
\begin{equation*}
\begin{aligned}
 \fint_{B_S(x)}\sigma(B_r(y))\,dy
 &=\int_{\R^N}\frac{|B_S(x)\cap B_r(z)|}{|B_S(x)|}\,d\sigma(z)
 \\&\le\int_{\R^N}
       \min\left\{1,\frac{|B_r(z)|}{|B_S(x)|}\right\}\,d\sigma(z)
 =\min\{1,(r/S)^N\}\,\sigma(\R^N).
\end{aligned}
\end{equation*}
Then, the Jensen's inequality gives that, for $0<L\le\infty$,
\begin{equation}\label{eq:average-trunc}
\begin{aligned}
 \fint_{B_S(x)}\V^L[\sigma](y)\,dy
 &\le\int_0^L r^{2-N/k}
   \left(\fint_{B_S(x)}\sigma(B_r(y))\,dy\right)^{1/k}\frac{dr}{r}\\
 &\le\sigma(\R^N)^{1/k}\int_0^L r^{2-N/k}
          \min\{1,(r/S)^{N/k}\}\frac{dr}{r}.
\end{aligned}
\end{equation}
For $L=\infty$, since $N>2k$, this gives
\begin{equation}\label{eq:exist-average}
\begin{aligned}
 \fint_{B_S(x)}\V[\sigma](y)\,dy
 &\le\sigma(\R^N)^{1/k}
    \left(S^{-N/k}\int_0^S r\,dr+
                 \int_S^\infty r^{1-N/k}\,dr\right)\\
 &=\left(\frac12+\frac1{N/k-2}\right)
       S^{2-N/k}\sigma(\R^N)^{1/k}.
\end{aligned}
\end{equation}
By \eqref{notation in thm1}, $\lambda_g$ is concentrated on $B_m$. Since $B_{1/m}(z)\subset B_{4m}$ for every $z\in B_m$, we have
\begin{equation*}
\begin{aligned}
 \sigma_g(\R^N)
 =\lambda_g(B_m)+4^N|B_1|\mu(B_m)
 \le (1+4^N|B_1|)\lambda_g(B_m),
\end{aligned}
\end{equation*}
Therefore, taking $\sigma=\sigma_g$ and $S=m$ in \eqref{eq:exist-average},
and using $|w_g|\le C\V[\sigma_g]$ in $B_{4m}$, we obtain,
for every $x\in B_m$,
\begin{equation}
\begin{aligned}\label{eq:average-lambda}
 \frac1m\fint_{B_m(x)}|w_g(y)|\,dy
 &\le \frac{C}{m}\fint_{B_m(x)}\V[\sigma_g](y)\,dy
 \le C m^{1-N/k}\sigma_g(\R^N)^{1/k}\\
 &\le C m^{1-N/k}\lambda_g(B_m)^{1/k}
 \le C\int_{2m}^{4m}
       \bigl(r^{k-N}\lambda_g(B_r(x))\bigr)^{1/k}\frac{dr}{r}\\
 &\le C\W[\lambda_g](x).
\end{aligned}
\end{equation}
The fourth inequality follows from $B_m\subset B_r(x)$ and
$r^{k-N}\ge(4m)^{k-N}$ for $2m\le r\le4m$. Since $g\in\mathcal K_m$, we have $|g|\le A\W[\mu]$ almost
everywhere in $B_m$. By the definition of $\lambda_g$ and
\eqref{eq:closure}, we obtain, almost everywhere in $B_m$,
\begin{equation}
\begin{aligned}\label{lambda potential estimate}
 \W[\lambda_g](x)
 &\le \W[\mu](x)+\W[\mathbf1_{B_m}|g|^q\,dx](x)\\
 &\le \W[\mu](x)+A^{q/k}\W[(\W[\mu])^q\,dx](x)\\
 &\le \bigl(1+CA^{q/k}D^{(q-k)/k^2}\bigr)\W[\mu](x).
\end{aligned}
\end{equation}
Therefore,
\begin{equation}\label{wg average estimate}
\frac1m\fint_{B_m(x)}|w_g(y)|\,dy\le C\bigl(1+CA^{q/k}D^{(q-k)/k^2}\bigr)\W[\mu](x).
\end{equation}
To estimate the remaining term in \eqref{eq:lmz}, it suffices to show
that $\W^m[\sigma_g](x)\le C\W[\lambda_g](x)$ for $x\in B_m$. Since $\sigma_g=\mathbf1_{B_{4m}}(\rho_{1/m}*\lambda_g+b_m)\,dx$
and $B_r(x)\subset B_{4m}$ for $0<r\le m$, we have
\begin{equation*}
\begin{aligned}
 \W^m[\sigma_g](x)
 &\le \W^m[(\rho_{1/m}*\lambda_g)\,dx](x)
      +\int_0^m\bigl(r^{k-N}b_m|B_r(x)|\bigr)^{1/k}\frac{dr}{r}\\
 &=\W^m[(\rho_{1/m}*\lambda_g)\,dx](x)
      +|B_1|^{1/k}m^{1-N/k}\mu(B_m)^{1/k}\\
 &\le \W^m[(\rho_{1/m}*\lambda_g)\,dx](x)
      +|B_1|^{1/k}m^{1-N/k}\lambda_g(B_m)^{1/k}\\
 &\le \W^m[(\rho_{1/m}*\lambda_g)\,dx](x)+C\W[\lambda_g](x),
\end{aligned}
\end{equation*}
where the last inequality follows from \eqref{eq:average-lambda}.
We now estimate $\W^m[(\rho_{1/m}*\lambda_g)\,dx](x)$.

Since $\operatorname{supp}\rho_h\subset B_h$, for any nonnegative finite measure $\sigma$ on $\mathbb{R}^N$,
\begin{equation}\label{measure mollifi}
 \int_{B_r(x)}(\rho_h*\sigma)(y)\,dy
 =\int_{B_{r+h}(x)}\int_{B_r(x)}\rho_h(y-z)\,dy\,d\sigma(z)
 \le \min\{1,|B_r|\|\rho_h\|_\infty\}\,\sigma(B_{r+h}(x)).
\end{equation}
Notice that for $r<h$ we have $B_{r+h}(x)\subset B_{2h}(x)$, and for $r\ge h$
we have $B_{r+h}(x)\subseteq B_{2r}(x)$. Thus, for $0<2h<S\le\infty$,
\begin{equation}\label{eq:exist-mollify}
\begin{aligned}
 \W^S[(\rho_h*\sigma)\,dx](x)
 &\le C\bigl(h^{-N}\sigma(B_{2h}(x))\bigr)^{1/k}\int_0^h dr
     +C\int_{2h}^{2S}\bigl(r^{k-N}\sigma(B_r(x))\bigr)^{1/k}\frac{dr}{r}\\
 &\le C\left(\int_{2h}^{4h}+\int_{2h}^{2S}\right)
           \bigl(r^{k-N}\sigma(B_r(x))\bigr)^{1/k}\frac{dr}{r}\le C\W^{2S}[\sigma](x).
\end{aligned}
\end{equation}
The second inequality follows from $\sigma(B_r(x))\ge\sigma(B_{2h}(x))$
for $2h\le r\le4h$.
Since $m\ge2$, we may take $\sigma=\lambda_g$, $h=1/m$, and $S=m$
in \eqref{eq:exist-mollify}. Thus, for $x\in B_m$,
\begin{equation}
\begin{aligned}\label{sigmag potential estimate}
 \W^m[\sigma_g](x)
 &\le \W^m[(\rho_{1/m}*\lambda_g)\,dx](x)+C\W[\lambda_g](x)\\
 &\le C\W^{2m}[\lambda_g](x)+C\W[\lambda_g](x)\\
 &\le C\W[\lambda_g](x).
\end{aligned}
\end{equation}
Substituting \eqref{sigmag potential estimate}, \eqref{lambda potential estimate} and \eqref{wg average estimate} into
\eqref{eq:lmz}, we obtain that

\begin{equation*}
\begin{aligned}
 |T_mg(x)|=|Dw_g(x)|\le C_0\bigl(1+A^{q/k}D^{(q-k)/k^2}\bigr)\W[\mu](x),
\end{aligned}
\end{equation*}
 almost everywhere in $B_m$. Here $C_0$ depends only on $N,k,q$ and the fixed mollifier. Choose $A=2C_0$ and $\varepsilon_0=A^{-kq/(q-k)}$.
For $D\le\varepsilon_0$ and $g\in\mathcal K_m$, we have,
almost everywhere in $B_m$,
\begin{equation*}
\begin{aligned}
 |T_mg|
 &\le C_0\bigl(1+A^{q/k}D^{(q-k)/k^2}\bigr)\W[\mu]\\
 &\le C_0\bigl(1+A^{q/k}\varepsilon_0^{(q-k)/k^2}\bigr)\W[\mu]
 =2C_0\W[\mu]
 =A\W[\mu].
\end{aligned}
\end{equation*}
Hence $T_m(\mathcal K_m)\subseteq\mathcal K_m$ for every such $m$.

We next show that $T_m(\mathcal K_m)$ is relatively compact in
$L^q(B_m;\R^N)$. To prove that, we want to show that for any $g\in\mathcal K_m$, $T_m g\in BV(B_m;\R^N)$.

Since $w_g$ is $k$-admissible and $k\ge2$, we have
$\Delta w_g\ge0$ and $\sigma_2(D^2w_g)\ge0$, so that
\begin{equation*}
 |D^2w_g|
 =\bigl((\Delta w_g)^2-2\sigma_2(D^2w_g)\bigr)^{1/2}
 \le\Delta w_g.
\end{equation*}
Fix $0\le\zeta\le1$ in $C_c^\infty(B_{2m})$ with $\zeta=1$ on $B_m$.
Since $T_mg=Dw_g$ in $B_m$ and $|Dw_g|\le A\W[\mu]$ almost everywhere
in $B_m$, we obtain
\begin{equation}
\begin{aligned}\label{eq:exist-bv}
 \|T_mg\|_{BV(B_m;\R^N)}
 &=\int_{B_m}|Dw_g|\,dx+\int_{B_m}|D^2w_g|\,dx\\
 &\le A\int_{B_m}\W[\mu]\,dx
      +\int_{B_{2m}}\zeta\Delta w_g\,dx\\
 &=A\int_{B_m}\W[\mu]\,dx
      +\int_{B_{2m}}w_g\Delta\zeta\,dx\\
 &\le A|B_m|^{1-1/q}\|\W[\mu]\|_{L^q(\R^N)}
      +\|\Delta\zeta\|_\infty |B_{2m}|
       \fint_{B_{2m}}|w_g|\,dx.
\end{aligned}
\end{equation}
So it remains to estimate the average of $|w_g|$. For every $B_S(x)\subset B_{4m}$, because $g\in\mathcal K_m$, \eqref{eq:energy}, \eqref{wg V estimate}
and \eqref{eq:exist-average} shows that
\begin{equation}\label{eq:exist-mean-uniform}
\begin{aligned}
 \fint_{B_S(x)}|w_g(y)|\,dy
 &\le C\fint_{B_S(x)}\V[\sigma_g](y)\,dy
 \le C S^{2-N/k}\sigma_g(\R^N)^{1/k}\\
 &=C S^{2-N/k}
   \left((1+4^N|B_1|)\mu(B_m)+\int_{B_m}|g|^q\,dy\right)^{1/k}\\
 &\le C S^{2-N/k}
   \left(\mu(\R^N)+A^q\int_{\R^N}(\W[\mu])^q\,dy\right)^{1/k}\\
 &\le C S^{2-N/k}
   \bigl(1+A^qD^{(q-k)/k}\bigr)^{1/k}\mu(\R^N)^{1/k}\\
 &\le C S^{2-N/k}\mu(\R^N)^{1/k},
\end{aligned}
\end{equation}
where the last inequality is obtained by $A^qD^{(q-k)/k}=(A^{q/k}D^{(q-k)/k^2})^k\le1$. Therefore, 
\begin{equation}
\|T_mg\|_{BV(B_m;\R^N)}\le C A|B_m|^{1-1/q}D^{(q-k)/(kq)}\mu(\R^N)^{1/q}+C_m\mu(\R^N)^{1/k}.
\end{equation}

Let $g_j\in\mathcal K_m$. By the compact embedding of $BV(B_m;\R^N)$ into $L^1(B_m;\R^N)$, after passing to a subsequence, there exists a $v\in L^1(B_m;\R^N)$
such that $T_mg_j\to v$ almost everywhere in $B_m$. Since
\begin{equation*}
 |T_mg_j-v|^q\le (2A\W[\mu])^q\in L^1(B_m),
\end{equation*}
dominated convergence theorem gives $T_mg_j\to v$ in $L^q(B_m;\R^N)$. Thus $T_m(\mathcal K_m)$ is relatively compact in $L^q(B_m;\R^N)$.

We next prove that $T_m$ is continuous. Let $g_j,g\in\mathcal K_m$ and suppose that $g_j\to g$ in $L^q(B_m;\R^N)$. By the definition of $f_g$ and H\"older's inequality,
\begin{equation*}
\begin{aligned}
 \|f_{g_j}-f_g\|_{L^\infty(B_{4m})}
 &\le \|\rho_{1/m}\|_\infty
      \bigl\||g_j|^q-|g|^q\bigr\|_{L^1(B_m)}\\
 &\le q\|\rho_{1/m}\|_\infty
      \bigl(\|g_j\|_{L^q(B_m)}+\|g\|_{L^q(B_m)}\bigr)^{q-1}
      \|g_j-g\|_{L^q(B_m)}
 \longrightarrow0,
\end{aligned}
\end{equation*}
as $j\to\infty$. To compare $w_{g_j}$ and $w_g$, choose $\tau_j\ge0$ such that
$\binom Nk\tau_j^k=\|f_{g_j}-f_g\|_{L^\infty(B_{4m})}$. Since $w_g$ is $k$-admissible and $\tau_j\ge0$, the function $w_g(x)+\frac{\tau_j}{2}(|x|^2-16m^2)$
is $k$-admissible. Moreover, since
$\sigma_i(D^2w_g)\ge0$ for $1\le i\le k$, we have
\begin{equation*}
\begin{aligned}
 \sigma_k\left(D^2\left(w_g+\frac{\tau_j}{2}(|x|^2-16m^2)\right)\right)
 &=\sum_{i=0}^k\binom{N-i}{k-i}\tau_j^{k-i}\sigma_i(D^2w_g)
 \ge \sigma_k(D^2w_g)+\binom Nk\tau_j^k\\
 &=f_g+\|f_{g_j}-f_g\|_{L^\infty(B_{4m})}
 \ge f_{g_j}
 =\sigma_k(D^2w_{g_j}),
\end{aligned}
\end{equation*}
where $\sigma_0=1$. On the other hand, since
$w_g=w_{g_j}=0$ on $\partial B_{4m}$, $w_g(x)+\frac{\tau_j}{2}(|x|^2-16m^2)=0$ on $\partial B_{4m}$. Therefore, the comparison principle gives
\begin{equation*}
 w_g(x)+\frac{\tau_j}{2}(|x|^2-16m^2)
 \le w_{g_j}(x)
 \quad\text{in }B_{4m}.
\end{equation*}
Interchanging $g$ and $g_j$, we obtain
\begin{equation*}
 -\frac{\tau_j}{2}(16m^2-|x|^2)
 \le w_{g_j}(x)-w_g(x)
 \le\frac{\tau_j}{2}(16m^2-|x|^2)
 \quad\text{in }B_{4m}.
\end{equation*}
Consequently, since
$\|f_{g_j}-f_g\|_{L^\infty(B_{4m})}\to0$, we have
\begin{equation}\label{eq:uniform-comparison}
 \|w_{g_j}-w_g\|_{L^\infty(B_{4m})}
 \le8m^2\tau_j\longrightarrow0.
\end{equation}

Let $\{T_mg_{j_\ell}\}$ be any subsequence converging to $v$
in $L^q(B_m;\R^N)$. For every $\phi\in C_c^\infty(B_m;\R^N)$,
integration by parts and \eqref{eq:uniform-comparison} give
\begin{equation*}
\begin{aligned}
 \int_{B_m}v\cdot\phi\,dx
 &=\lim_{\ell\to\infty}\int_{B_m}Dw_{g_{j_\ell}}\cdot\phi\,dx
 =-\lim_{\ell\to\infty}
     \int_{B_m}w_{g_{j_\ell}}\operatorname{div}\phi\,dx\\
 &=-\int_{B_m}w_g\operatorname{div}\phi\,dx
 =\int_{B_m}Dw_g\cdot\phi\,dx.
\end{aligned}
\end{equation*}
Thus $v=Dw_g=T_mg$ almost everywhere in $B_m$. Since
$T_m(\mathcal K_m)$ is relatively compact in $L^q(B_m;\R^N)$,
it follows that $T_mg_j\to T_mg$ in $L^q(B_m;\R^N)$.

By Schauder's fixed point theorem, there exists $g_m\in\mathcal K_m$ such that $T_mg_m=g_m$. Set $w_m=w_{g_m}$ and extend $g_m$ by zero outside $B_m$. By the definition of $T_m$ and $\mathcal K_m$, we have
\begin{equation*}
 g_m=Dw_m\quad\text{a.e. in }B_m,
 \quad
 |g_m|\le A\W[\mu]\quad\text{a.e. in }\R^N.
\end{equation*}
Moreover, by the definition of $w_g$ and $f_g$,
\begin{equation*}
 \sigma_k(D^2w_m)
 =\rho_{1/m}*
 \left(\mathbf1_{B_m}\mu
 +\mathbf1_{B_m}|Dw_m|^q\,dx\right)+b_m
 \quad\text{in }B_{4m}.
\end{equation*}

We now let $m\to\infty$. For every fixed $S>0$,
arguing as in \eqref{eq:exist-bv} with a cutoff equal to one
on $B_S$ and supported in $B_{2S}$, and using
\eqref{eq:exist-mean-uniform}, \eqref{eq:energy},
and $|g_m|\le A\W[\mu]$, we obtain
\begin{equation*}
 \sup_{m>2S}
 \left(
 \|w_m\|_{W^{1,1}(B_S)}
 +\|Dw_m\|_{BV(B_S;\R^N)}
 \right)<\infty.
\end{equation*}
By $BV$ compactness and a diagonal argument,
we may choose a common subsequence, still denoted by $\{w_m\}$ and $\{g_m\}$, such that
$w_m\to w$ and $g_m\to Dw$ in $L^1_{\rm loc}(\R^N)$ and almost everywhere in $\R^N$.  Since 
$|g_m|\le A\W[\mu]$, by \eqref{eq:energy} and
the dominated convergence theorem, we obtain
\begin{equation*}
\begin{aligned}
 |Dw|&\le A\W[\mu]\quad\text{a.e. in }\R^N,\\
 \|g_m-Dw\|_{L^q(\R^N)}
 &+\bigl\||g_m|^q-|Dw|^q\bigr\|_{L^1(\R^N)}
 \longrightarrow0.
\end{aligned}
\end{equation*}
Moreover, by Poincar\'e's inequality and
\eqref{eq:exist-mean-uniform}, for every fixed $S>0$,
\begin{equation*}
\begin{aligned}
 \sup_{m>2S}\|w_m\|_{L^q(B_S)}
 &\le \sup_{m>2S}\left(\|w_m-(w_m)_{B_S}\|_{L^q(B_S)}
   +|B_S|^{1/q}|(w_m)_{B_S}|\right)\\
 &\le C_S\sup_{m>2S}
 \left(\|w_m\|_{L^1(B_S)}
       +\|Dw_m\|_{L^q(B_S)}\right)\\
 &\le C_S\left(
 \mu(\R^N)^{1/k}
 +A\|\W[\mu]\|_{L^q(\R^N)}\right)<\infty.
\end{aligned}
\end{equation*}
Since $w_m\to w$ almost everywhere, Fatou's lemma gives
$w\in L^q_{\rm loc}(\R^N)$. Together with $Dw\in L^q(\R^N)$, this yields $w\in W^{1,q}_{\rm loc}(\R^N)$. Since $w_m$ is $k$-convex and $w_m\to w$ in
$L^1_{\rm loc}(\R^N)$, \cite[Lemma~2.2]{TW}
gives a $k$-convex representative of $w$.
Moreover, since $w_m\le0$, the submean inequality
shows that $w\le0$ everywhere.

By the weak continuity of Hessian measures
\cite[Theorem~1.1(ii)]{TW}, for every
$\phi\in C_c(\R^N)$, we obtain
\begin{equation*}
\begin{aligned}
 \int_{\R^N}\phi\,dF_k[w]
 &=\lim_{m\to\infty}
   \int_{B_{4m}}\phi f_{g_m}\,dx\\
 &=\lim_{m\to\infty}
   \left(
   \int_{B_m}(\rho_{1/m}*\phi)\,d\mu
   +\int_{\R^N}(\rho_{1/m}*\phi)|g_m|^q\,dx
   +b_m\int_{\R^N}\phi\,dx
   \right)\\
 &=\int_{\R^N}\phi\,d\mu
   +\int_{\R^N}\phi|Dw|^q\,dx,
\end{aligned}
\end{equation*}
where the last equality follows from
$\rho_{1/m}*\phi\to\phi$ uniformly,
$|g_m|^q\to|Dw|^q$ in $L^1(\R^N)$,
and $0\le b_m\le m^{-N}\mu(\R^N)\to0$.
Therefore,
\begin{equation*}
 F_k[w]=|Dw|^q\,dx+\mu
 \quad\text{in }\R^N.
\end{equation*}
Set $u=-w$. We have $F_k[-u]=|Du|^q\,dx+\mu$
and \eqref{eq:global-gradient}. By
\eqref{eq:exist-mean-uniform}, for every $S>0$,
\begin{equation*}
\begin{aligned}
 0\le\inf_{\R^N}u
 &\le\fint_{B_S}u\,dx
 =\lim_{m\to\infty}\fint_{B_S}|w_m|\,dx
 \le CS^{2-N/k}\mu(\R^N)^{1/k}
 \longrightarrow0,
\end{aligned}
\end{equation*}
as $S\to\infty$. Hence $\inf_{\R^N}u=0$.

This proves the assertion for $t=1$ when
$D\le\varepsilon_0$.
For arbitrary finite $D>0$, take
$0<t\le\varepsilon_0/D$ and apply the preceding
result to $t\mu$.
Conversely, Proposition~\ref{prop:hessian-capacity}
gives
\begin{equation*}
 t\mu(K)\le F_k[-u](K)
 \le C(N,k,q)\Capacity_{I_k,q/(q-k)}(K)
\end{equation*}
for every compact $K\subset\R^N$.
This proves \eqref{eq:riesz-condition}
and completes the proof.
\end{proof}

\begin{proof}[Proof of Corollary~\ref{cor:locally-finite}]
Let $s=q/(q-k)$ and write $B_R=B_R(0)$.
Suppose first that a solution exists.
By Proposition~\ref{prop:hessian-capacity},
\eqref{eq:riesz-condition} holds with $D=C/t$.
Choose $x\in\R^N$ with $u(x)<\infty$ and let $R=1+|x|$.
Since $B_r(0)\subset B_{2r}(x)$ for $r\ge R$,
\cite[Corollary~7.3]{PV} gives
\begin{equation*}
\begin{aligned}
 \int_R^\infty
 \bigl(r^{2k-N}\mu(B_r(0))\bigr)^{1/k}\frac{dr}{r}
 &\le C\int_{2R}^\infty
 \bigl(r^{2k-N}\mu(B_r(x))\bigr)^{1/k}\frac{dr}{r}\\
 &\le C\V[\mu](x)
 \le Ct^{-1/k}u(x)<\infty.
\end{aligned}
\end{equation*}
Thus \eqref{eq:tail-locally-finite} holds.

If $q<2k$, then $s>2$ and \eqref{eq:riesz-condition} gives
$\mu(B_r)\le CDr^{N-ks}$. Hence
\begin{equation*}
\begin{aligned}
 \int_1^\infty
 \bigl(r^{2k-N}\mu(B_r)\bigr)^{1/k}\frac{dr}{r}
 \le CD^{1/k}\int_1^\infty r^{2-s}\frac{dr}{r}
 =\frac{CD^{1/k}}{s-2}<\infty.
\end{aligned}
\end{equation*}

Conversely, suppose that \eqref{eq:riesz-condition} and \eqref{eq:tail-locally-finite} hold. If $\mu=0$, the conclusion follows by taking $u\equiv0$. We may therefore assume that $\mu\ne0$. We first consider the case $D\le\varepsilon_0$.

Set $\mu_j=\mathbf1_{B_j}\mu$. Since $\mu$ is locally finite,
each $\mu_j$ is finite. Moreover, for every compact set
$K\subset\R^N$, we have
\begin{equation*}
 \mu_j(K)\le\mu(K)\le D\Capacity_{I_k,s}(K).
\end{equation*}
Since $D\le\varepsilon_0$, Theorem~\ref{thm:exist} gives a nonnegative
solution $u_j\in W^{1,q}_{\rm loc}(\R^N)$ with
$-u_j\in\Phi_k(\R^N)$ and $\inf_{\R^N}u_j=0$ such that
\begin{equation*}
 F_k[-u_j]=|Du_j|^q\,dx+\mu_j \quad\text{ in }\R^N.
\end{equation*}
Moreover, the gradient of this solution satisfies that
\begin{equation*}
 |Du_j|\le C\W[\mu_j]\le C\W[\mu]
 \quad\text{a.e. in }\R^N.
\end{equation*}
We now derive an upper bound for $u_j$ independent of $j$. By \cite[Theorem~2.1(a),(b)]{Mixed}, the capacity condition \eqref{eq:riesz-condition} implies that $\W[\mu]\in L^q_{\rm loc}(\R^N)$. Then
\begin{equation}\label{eq:locally-finite-upper}
\begin{aligned}
 u_j
 &\le C\V[|Du_j|^q\,dx+\mu_j]
 \le C\V[\mu]+C\V[(\W[\mu])^q\,dx]\\
 &=C\V[\mu]
   +CD^{q/k^2}\V[(\W[\mu/D])^q\,dx]\\
 &\le C\V[\mu]+CD^{q/k^2}\V[\mu/D]\\
 &=C\bigl(1+D^{(q-k)/k^2}\bigr)\V[\mu]\\
 &\le C\V[\mu]
 \quad\text{a.e. in }\R^N.
\end{aligned}
\end{equation}
The first inequality follows from
\cite[Corollary~7.3]{PV}, the third inequality follows by applying \cite[Theorem~2.1 and Proposition~2.4(i)]{Mixed} to $\mu/D$, and the last inequality follows from $D\le\varepsilon_0$.

To obtain the bounds independent of $j$ on each ball,
we estimate the averages of $\V[\mu]$.
For $x\in\R^N$ and $r,R>0$, direct computation gives
\begin{equation*}
\begin{aligned}
 \fint_{B_R(x)}\mu(B_r(y))\,dy
 &=\int_{\R^N}
 \frac{|B_R(x)\cap B_r(z)|}{|B_R(x)|}\,d\mu(z)\le\min\{1,(r/R)^N\}\mu(B_{R+r}(x)).
\end{aligned}
\end{equation*}
By Jensen's inequality,
\begin{equation}\label{eq:locally-finite-average}
\begin{aligned}
 \fint_{B_R(x)}\V[\mu](y)\,dy
 &\le\int_0^\infty r^{2-N/k}
 \left(\fint_{B_R(x)}\mu(B_r(y))\,dy\right)^{1/k}
 \frac{dr}{r}\\
 &\le CR^{2-N/k}\mu(B_{2R}(x))^{1/k}
 +C\int_{2R}^\infty
 \bigl(r^{2k-N}\mu(B_r(x))\bigr)^{1/k}\frac{dr}{r}\\
 &\le C\int_{2R}^\infty
 \bigl(r^{2k-N}\mu(B_r(x))\bigr)^{1/k}\frac{dr}{r}.
\end{aligned}
\end{equation}
Taking $x=0$ and using \eqref{eq:tail-locally-finite},
we obtain $\V[\mu]\in L^1_{\rm loc}(\R^N)$ and
$
 \fint_{B_R}\V[\mu]\,dx\longrightarrow0,
$
as $R\to\infty$.

We next obtain a uniform local $BV$ bound for $Du_j$. The standard mollifications of $-u_j$ are smooth and
$k$-convex by \cite[Lemma~2.3]{TW}. We apply the integration by parts argument leading to \eqref{eq:exist-bv} to these mollifications. Passing to the limit by lower semicontinuity and using \eqref{eq:locally-finite-upper} and $|Du_j|\le C\W[\mu]$, we obtain, for every $R>0$,
\begin{equation*}
 \sup_j\left(
 \|u_j\|_{W^{1,1}(B_R)}
 +\|Du_j\|_{BV(B_R;\R^N)}
 \right)
 \le C_R\int_{B_{2R}}
       \bigl(\V[\mu]+\W[\mu]\bigr)\,dx<\infty,
\end{equation*}
where $C_R$ is independent of $j$.

Arguing as the proof of Theorem~\ref{thm:exist}, by $BV$ compactness,
Poincar\'e's inequality, and dominated convergence with
$\W[\mu]\in L^q_{\rm loc}(\R^N)$, we may choose a subsequence,
still denoted by $\{u_j\}$, and
$u\in W^{1,q}_{\rm loc}(\R^N)$ such that
$u_j\to u$ and $Du_j\to Du$ almost everywhere, and
\begin{equation*}
\begin{aligned}
 u_j&\to u
       &&\text{in }L^1_{\rm loc}(\R^N),\\
 Du_j&\to Du
       &&\text{in }L^q_{\rm loc}(\R^N;\R^N),\\
 |Du_j|^q&\to|Du|^q
       &&\text{in }L^1_{\rm loc}(\R^N).
\end{aligned}
\end{equation*}
Moreover, we can also get that $0\le u\le C\V[\mu]$ and
$|Du|\le C\W[\mu]$ almost everywhere.
By \cite[Lemma~2.2]{TW}, we take $-u$ to be its
$k$-convex representative; the mean value inequality
gives $u\ge0$ everywhere.
Since $\mu_j=\mathbf1_{B_j}\mu\rightharpoonup\mu$ locally
as measures, \cite[Theorem~1.1(ii)]{TW} gives
\begin{equation*}
 F_k[-u]
 =\lim_{j\to\infty}F_k[-u_j]
 =\lim_{j\to\infty}\bigl(|Du_j|^q\,dx+\mu_j\bigr)
 =|Du|^q\,dx+\mu
\end{equation*}
in the sense of measures.

By \eqref{eq:locally-finite-average} with $x=0$
and \eqref{eq:tail-locally-finite},
\begin{equation*}
 0\le\inf_{\R^N}u
 \le\fint_{B_R}u\,dx
 \le C\fint_{B_R}\V[\mu]\,dx
 \longrightarrow0
 \quad\text{as }R\to\infty.
\end{equation*}
Thus $\inf_{\R^N}u=0$, and \cite[Corollary~7.3]{PV}
gives $u(x)\ge c\V[\mu](x)$ for every $x\in\R^N$.
For the upper bound, since $u$ is superharmonic and
$u\le C\V[\mu]$ almost everywhere,
\eqref{eq:locally-finite-average} shows that for every $x\in\R^N$,
\begin{equation*}
\begin{aligned}
 u(x)
=\lim_{r\downarrow0}\fint_{B_r(x)}u(y)\,dy
 \le C\lim_{r\downarrow0}\int_{2r}^\infty
       \bigl(t^{2k-N}\mu(B_t(x))\bigr)^{1/k}\frac{dt}{t}
 =C\V[\mu](x).
\end{aligned}
\end{equation*}

For arbitrary $D>0$, apply the preceding argument to $t\mu$
with $0<t\le\varepsilon_0/D$.
The constant in \eqref{eq:riesz-condition} becomes $tD$,
and the integral in \eqref{eq:tail-locally-finite}
is multiplied by $t^{1/k}$.
This completes the proof.
\end{proof}
\begin{remark}
For $q\ge2k$, \eqref{eq:riesz-condition} does not imply
\eqref{eq:tail-locally-finite}.
In fact, let $s=q/(q-k)$ and $d\mu(x)=|x|^{-ks}\,dx$.
By the Sobolev inequality for Riesz potentials,
\begin{equation*}
 \mu(K)\le C|K|^{1-ks/N}
 \le C\Capacity_{I_k,s}(K)
\end{equation*}
for every compact $K\subset\R^N$. However,
$\mu(B_r)=Cr^{N-ks}$ and $s\le2$, so
\begin{equation*}
 \int_1^\infty
 \left(\frac{\mu(B_r)}{r^{N-2k}}\right)^{1/k}
 \frac{dr}{r}
 =C\int_1^\infty r^{1-s}\,dr=\infty.
\end{equation*}
\end{remark}

\section{Fully Nonlinear Dirichlet problem}\label{sec:domains}
\begin{proof}[Proof of Theorem~\ref{thm:ball}]
We follow the proof of Theorem~\ref{thm:exist} and give
the modifications needed for the bounded domain.

We first take $t=1$ and suppose that $D$ is sufficiently small.
Write $s=q/(q-k)$, $L=8\operatorname{diam}\Omega$,
$M=\mu(\Omega)$ and $d(x)=\operatorname{dist}(x,\partial\Omega)$.
If $M=0$, take $u=0$. Assume $M>0$.
Since $\operatorname{supp}\mu\subset E$, condition
\eqref{eq:bessel-condition} extends to every compact
$K\subset\R^N$, and $M\le D\Capacity_{G_k,s}(E)$.
For $x_0\in\Omega$, we have $E\subset B_{L/8}(x_0)$.
Lemma~\ref{lem:potential} therefore gives \eqref{eq:energy}
and \eqref{eq:closure}.

Choose $0<\delta<L/12$ such that $d$ is smooth in
$\{0<d<4\delta\}\cap\Omega$ and
$6\delta<\operatorname{dist}(E,\partial\Omega)$.
Then
\begin{equation}\label{eq:collar-potential}
\begin{aligned}
 \W^L[\mu]&\ge cM^{1/k}
     &&\text{in }\Omega,\\
 \W^L[\mu]&\le CM^{1/k}
     &&\text{on }\{d<4\delta\}\cap\Omega.
\end{aligned}
\end{equation}
Throughout the proof, $C$ denotes a positive constant
depending only on $N,k,q,\Omega,E$.

For $A\ge1$ and $0<h<\min\{1,\delta/8\}$, set
\begin{equation*}
 \mathcal K=\{g\in L^q(\Omega;\R^N):
       |g|\le A\W^L[\mu]\text{ a.e. in }\Omega\}.
\end{equation*}
For $g\in\mathcal K$, set
$\lambda_g=\mu+\mathbf1_\Omega|g|^q\,dx$ and
$f_g=\rho_h*\lambda_g+hM$, where $\rho_h$ is the mollifier
fixed in the proof of Theorem~\ref{thm:exist}.
By \cite[Theorem~1]{CNS}, there exists a smooth
$k$-admissible solution $w_g$ of
\begin{equation*}
 \sigma_k(D^2w_g)=f_g\quad\text{in }\Omega,
 \quad w_g=0\quad\text{on }\partial\Omega.
\end{equation*}
Set $T_hg=Dw_g$.

Suppose $A^qD^{(q-k)/k}\le1$.
For $d(x)<3\delta$, we have $(\rho_h*\mu)(x)=0$ and
$B_h(x)\cap\Omega\subset\{d<4\delta\}\cap\Omega$.
By \eqref{measure mollifi}, \eqref{eq:energy}
and \eqref{eq:collar-potential},
\begin{equation}\label{eq:domain-mass-collar}
\begin{aligned}
 \int_\Omega f_g\,dx
 &\le M+CA^qD^{(q-k)/k}M+hM|\Omega|
 \le CM,\\
 \sup_{\{d<3\delta\}\cap\Omega}f_g
 &\le CA^qM^{q/k}+hM
 \le CM.
\end{aligned}
\end{equation}
The last inequality in the second line uses
$M\le D\Capacity_{G_k,s}(E)$.

For $x\in\Omega$ with $d(x)\ge\delta/2$, apply
\eqref{eq:lmz} with radius $\delta/4$.
Arguing as in the proof of Theorem~\ref{thm:exist},
with \eqref{eq:average-trunc} in place of
\eqref{eq:exist-average}, using the first inequality of $\eqref{eq:collar-potential}$ and \eqref{eq:domain-mass-collar}, we obtain
\begin{equation}\label{eq:domain-interior}
 |Dw_g(x)|
 \le C_0\bigl(1+A^{q/k}D^{(q-k)/k^2}\bigr)
       \W^L[\mu](x)
\end{equation}
It remains to estimate $Dw_g$ near the boundary.
By \cite[Theorem~7.4]{PV} and
\eqref{eq:domain-mass-collar}, for $d(x)=\delta$,
\begin{equation*}
\begin{aligned}
 |w_g(x)|
 &\le C\V^{L/4}[\mathbf1_\Omega f_g\,dx](x)\\
 &=C\left(\int_0^{\delta/2}
       +\int_{\delta/2}^{L/4}\right)
   \left(t^{2k-N}\int_{B_t(x)\cap\Omega}f_g(y)\,dy
   \right)^{1/k}\frac{dt}{t}\\
 &\le CM^{1/k}\left(
       \int_0^{\delta/2}t^2\frac{dt}{t}
       +\int_{\delta/2}^{L/4}
          t^{2-N/k}\frac{dt}{t}\right)\\
 &\le CM^{1/k}.
\end{aligned}
\end{equation*}
Thus
\begin{equation}
\sup_{\{d=\delta\}}|w_g|\le CM^{1/k}. 
\end{equation}
We use this estimate to construct a barrier for $w_g$
in $\{0<d<\delta\}\cap\Omega$.
By \cite[Theorem~1]{CNS}, there exists a smooth
$k$-admissible solution $\psi$ satisfying
\begin{equation*}
 \sigma_k(D^2\psi)=1\quad\text{in }\Omega,
 \quad \psi=0\quad\text{on }\partial\Omega.
\end{equation*}
By the strong maximum principle, $\psi<0$ in $\Omega$.
Since $\psi\in C^\infty(\overline\Omega)$, we have
$\min_{\{d=\delta\}}(-\psi)>0$ and
$-\psi\le\|D\psi\|_{L^\infty(\Omega)}d$.
Thus, for a sufficiently large fixed $B>0$,
\begin{equation*}
 \sigma_k\bigl(D^2(BM^{1/k}\psi)\bigr)
 =B^kM\ge f_g
 \quad\text{in }\{0<d<\delta\}\cap\Omega,
\end{equation*}
and
\begin{equation*}
 BM^{1/k}\psi\le w_g
 \quad\text{on }\partial\bigl(\{0<d<\delta\}\cap\Omega\bigr).
\end{equation*}
Then, for every $x\in\Omega$ with $0<d(x)<\delta$, the comparison principle gives
\begin{equation}\label{eq:domain-barrier}
 0\le-w_g(x)
 \le-BM^{1/k}\psi(x)
 \le CM^{1/k}d(x)
.
\end{equation}
For $d(x)<\delta/2$, take $r=d(x)/2$.
Applying \eqref{eq:lmz} and using
\eqref{eq:domain-barrier} and \eqref{eq:domain-mass-collar},
we obtain
\begin{equation}\label{eq:domain-boundary-gradient}
\begin{aligned}
 |Dw_g(x)|
 &\le \frac{CM^{1/k}}r\fint_{B_r(x)}d(y)\,dy
       +CM^{1/k}r\\
 &\le CM^{1/k}
 \le C\W^L[\mu](x),
\end{aligned}
\end{equation}
the last inequality follows from the first inequality in
\eqref{eq:collar-potential}. Together with \eqref{eq:domain-interior}, this yields
\begin{equation*}
 |T_hg|
 \le C_1\bigl(1+A^{q/k}D^{(q-k)/k^2}\bigr)
       \W^L[\mu]
 \quad\text{a.e. in }\Omega.
\end{equation*}
Choose $A=\max\{1,2C_1\}$ and $\varepsilon_0>0$ such that
$A^q\varepsilon_0^{(q-k)/k}\le1$.
Then $T_h(\mathcal K)\subset\mathcal K$ for
$D\le\varepsilon_0$, as in Theorem~\ref{thm:exist}.

For every fixed $U\Subset\Omega$, arguing as in
\eqref{eq:exist-bv} with a cutoff equal to one on $U$
and supported in $\Omega$, and using Poincar\'e's
inequality, \eqref{eq:energy}, and
$|Dw_g|\le A\W^L[\mu]$, we obtain
\begin{equation*}
 \sup_{g\in\mathcal K}
 \left(
 \|w_g\|_{W^{1,q}_0(\Omega)}
 +\|Dw_g\|_{BV(U;\R^N)}
 \right)<\infty,
\end{equation*}
uniformly in $h$. As in the proof of Theorem~\ref{thm:exist}, the compact embedding of $BV(U;\R^N)$ into
$L^1(U;\R^N)$ for every smooth $U\Subset\Omega$, together with a diagonal argument, shows that
every sequence in $\bigcup_{0<h<h_0}T_h(\mathcal K)$
has an almost everywhere convergent subsequence
in $\Omega$. Since
$|T_hg|\le A\W^L[\mu]\in L^q(\Omega)$,
the dominated convergence theorem gives
the relative compactness of
$\bigcup_{0<h<h_0}T_h(\mathcal K)$
in $L^q(\Omega;\R^N)$.

Following the proof of Theorem~\ref{thm:exist}. Fix $h>0$ and let $g_j\to g$ in $L^q(\Omega;\R^N)$, with
$g_j,g\in\mathcal K$. Then $f_{g_j}\to f_g$ uniformly in $\Omega$.
Choose $B_R(x_0)$ containing $\overline\Omega$ and set
$Q(x)=(|x-x_0|^2-R^2)/2$. As in the proof of
\eqref{eq:uniform-comparison}, we compare $w_g+\tau_jQ$ with $w_{g_j}$ and then interchange $g$ and $g_j$, where
$\tau_j=(\|f_{g_j}-f_g\|_\infty/\binom Nk)^{1/k}$.
Since $Q\le0$ on $\partial\Omega$, the two comparisons show that
$\|w_{g_j}-w_g\|_\infty\le R^2\tau_j/2\to0$.
Thus every $L^q$ subsequential limit of $Dw_{g_j}$ equals $Dw_g$
in distributions. The relative compactness proved above yields
$T_hg_j\to T_hg$ in $L^q(\Omega;\R^N)$.
Schauder's theorem now gives a fixed point $g_h=T_hg_h$.

Set $w_h=w_{g_h}$ and let $h\to0$. By the relative compactness proved above, we may pass to a subsequence, still denoted by $\{w_h\}$, such that $Dw_h$ converges strongly in $L^q(\Omega;\R^N)$. Since $w_h\in W^{1,q}_0(\Omega)$, Poincar\'e's inequality implies that $w_h$ converges strongly to some $w$ in $W^{1,q}_0(\Omega)$
along the same subsequence. As in the proof of Theorem~\ref{thm:exist}, it follows from \cite[Lemma~2.2 and Theorem~1.1(ii)]{TW} and the strong convergence
of $Dw_h$ that $w$ has a nonpositive $k$-convex representative and
satisfies
\begin{equation*}
 F_k[w]=|Dw|^q\,dx+\mu
\end{equation*}
in $\Omega$, with $|Dw|\le A\W^L[\mu]$ almost everywhere in $\Omega$. Moreover, $w$ also satisfies the bounds in \eqref{eq:domain-barrier} almost everywhere. Thus $u=-w$ is locally Lipschitz near $\partial\Omega$.
For every $x\in\Omega$ with $d(x)<\delta/2$, we have
\begin{equation*}
 0\le u(x)\le CM^{1/k}d(x).
\end{equation*}
Consequently,
$\lim_{\Omega\ni x\to z}u(x)=0$
for every $z\in\partial\Omega$.

For $0<D<\infty$, we apply the preceding argument to $t\mu$ with $0<t\le\varepsilon_0/D$.
For the converse, Proposition~\ref{prop:hessian-capacity} shows that
\begin{equation*}
 t\mu(K)\le F_k[-u](K)
 \le C(N,k,q,\Omega,E)\Capacity_{G_k,q/(q-k)}(K),
\end{equation*}
where $ K\subset E$ is a compact set. Thus \eqref{eq:bessel-condition} holds with $D=C(N,k,q,\Omega,E)/t$.
\end{proof}
\medskip
\noindent\textbf{Acknowledgments.}\enspace
The first named author is supported by the Natural Science Foundation of Tianjin, No. 22JCJQJC00130 
\\[0.5em]\noindent\textbf{AI Disclosure.}\enspace
The authors used ChatGPT to improve the language and assist in checking
mathematical arguments. The authors take full responsibility for the
final manuscript.

\bibliographystyle{amsplain}

\bibliography{references}

\end{document}